\documentclass{article}

\usepackage{amsmath,amssymb,amsthm}

\newtheorem{theorem}{Theorem}
\newtheorem{corollary}{Corollary}
\newtheorem{lemma}{Lemma}

\theoremstyle{remark}
\newtheorem{remark}{Remark}

\theoremstyle{definition}
\newtheorem{definition}{Definition}

\begin{document}

\title{Generalizations of Gronwall-Bihari\\ 
Inequalities on Time Scales\thanks{This is a preprint of an article accepted (16/May/2008) for publication in the \emph{Journal of Difference Equations and Applications} [\copyright\ Taylor \& Francis]; \emph{J. Difference Equ. Appl.} is available online at {\tt http://www.informaworld.com}.}}

\author{Rui A. C. Ferreira\thanks{Supported by FCT through the PhD fellowship SFRH/BD/39816/2007.}\\
{\tt ruiacferreira@ua.pt}
\and Delfim F. M. Torres\thanks{Supported by FCT through the R\&D unit CEOC, cofinanced by the
EC fund FEDER/POCI 2010.}\\ 
{\tt delfim@ua.pt}}

\date{Department of Mathematics\\ 
University of Aveiro\\ 
3810-193 Aveiro, Portugal}

\maketitle


\begin{abstract}
We establish some nonlinear integral
inequalities for functions defined on a time scale. The
results extend some previous Gronwall and Bihari type inequalities on time scales. Some examples of time scales for which our results can be applied are provided. An application to the qualitative analysis of a nonlinear dynamic equation is discussed.

\bigskip

\noindent \textbf{Keywords:} Inequalities, Gronwall and Bihari inequalities, subadditive and submultiplicative functions, time scales.

\bigskip

\noindent \textbf{2000 Mathematics Subject Classification:} 26D15, 34A40, 39A12.
\end{abstract}


\section{Introduction}

The theory of time scales was introduced in 1988 at
Stefan Hilger's PhD thesis, with the primary goal to unify and extend the continuous and discrete analysis \cite{Hilger90,Hilger97}.
Since then, the theory has been growing up and applied to many
different fields of mathematics \cite{Agarwal,livro,advance}.
We refer the reader to \cite{livro} for all the
basic definitions and results on time scales necessary to this work (\textrm{e.g.}, delta differentiability, rd-continuity,
exponential function and its properties).

It is well known that inequalities play an important role in the
study of differential and difference equations \cite{inesurvey,HLP,RachidDelfim}.
Among many types of important
inequalities is Gronwall inequality and their nonlinear extensions, namely Bihari type ones \cite{Beesack,deo,og}. These and many other types of inequalities have
been derived for the more general setting of time scales \cite{inesurvey,Pach,Bihary,Deepak,AFT,Gronwall}.

Motivated by the recent paper
\cite{motivacao}, we establish here some new nonlinear integral
inequalities on time scales. Our inequalities differ from those
found in the literature by the introduction of new kind of nonlinearities. In \cite{inesurvey,livro,Bihary} the assumption
\begin{equation}
\label{ineq:gi}
u(t)\leq a(t)+\int_a^t f(s)u(s)\Delta s
\end{equation}
is considered with respect to Gronwall's inequality (see \textrm{e.g.}, \cite[Theorem~5.6]{inesurvey});
the assumption
\begin{equation}
\label{ineq:bi}
u(t)\leq a(t)+\int_a^t f(s)g(u(s))\Delta s
\end{equation}
with respect to Bihari's inequality (see \textrm{e.g.}, \cite[Theorem~5.8]{inesurvey}).
Here we relax both hypotheses \eqref{ineq:gi} and \eqref{ineq:bi}
by adding to their right-hand sides new nonnegative terms (\textrm{cf.} \eqref{eq0} and \eqref{eq8}, respectively).
In Section~\ref{sec:mainResults} we state and prove our results; in Section~\ref{aplic} one of them is used to estimate the solution of a nonlinear dynamic equation.
We employ the concepts of delta-derivative and delta-integral; analogous results can be easily obtained using the so called nabla-derivative and nabla-integral (for definitions, see \cite{advance}). To the best of our knowledge, the
results are new even for the discrete time case, when the time scale is chosen to be the set of integers.


\section{Main Results}
\label{sec:mainResults}

Throughout we use the notation
$\mathbb{R}^+_0=[0,\infty)$. Let $\mathbb{T}$ be a time scale.
For $a,b\in\mathbb{T}$ with $a<b$, we define the time scales interval by $$[a,b]_{\mathbb{T}}=\{t\in\mathbb{T}:a\leq t\leq b\}.$$ Lemma~\ref{lemimp} is a useful tool for the proofs of the next theorems.

\begin{lemma}
\label{lemimp} Let $a,b\in\mathbb{T}$, consider the time scales interval $[a,b]_{\mathbb{T}}$, and a delta differentiable function $r:[a,b]_{\mathbb{T}}\rightarrow(0,\infty)$ with
$r^\Delta(t)\geq 0$ on $[a,b]_{\mathbb{T}}^{\kappa}$. Define
\begin{equation*}
G(x)=\int_{x_0}^{x}\frac{ds}{g(s)},\ x>0,\ x_0>0 \, ,
\end{equation*}
where $g\in C(\mathbb{R}^+_0,\mathbb{R}^+_0)$ is positive and
nondecreasing on $(0,\infty)$. Then, for each
$t\in[a,b]_{\mathbb{T}}$ we have
\begin{equation*}
G(r(t))\leq
G(r(a))+\int_{a}^{t}\frac{r^\Delta(\tau)}{g(r(\tau))}\Delta\tau.
\end{equation*}
\end{lemma}

\begin{proof}
Since $g$ is positive and nondecreasing on $(0,\infty)$, we have,
successively, that
\begin{equation}
\label{seila8}
\begin{gathered}
r(t) \leq r(t)+h\mu(t)r^\Delta(t) \, ,\\
g(r(t)) \leq g(r(t)+h\mu(t)r^\Delta(t)) \, ,\\
\frac{1}{g(r(t)+h\mu(t)r^\Delta(t))} \leq\frac{1}{g(r(t))} \, ,\\
\int_0^1\frac{1}{g(r(t)+h\mu(t)r^\Delta(t))}dh \leq\int_0^1\frac{1}{g(r(t))}dh=\frac{1}{g(r(t))} \, , \\
\left\{\int_0^1\frac{1}{g(r(t)+h\mu(t)r^\Delta(t))}dh\right\}
r^\Delta(t) \leq\frac{r^\Delta(t)}{g(r(t))} \, ,
\end{gathered}
\end{equation}
for all $t\in[a,b]_{\mathbb{T}}^{\kappa}$ and $h\in[0,1]$.
By $\Delta$-integrating the last inequality in (\ref{seila8}) from $a$ to $t$ and
having in mind that the chain rule \cite[Theorem~1.90]{livro}
guarantees that
\begin{align}
(G\circ r)^\Delta(t)&=\left\{\int_0^1
G'(r(t)+h\mu(t)r^\Delta(t))dh\right\}r^\Delta(t)\nonumber\\
&=\left\{\int_0^1
\frac{1}{g(r(t)+h\mu(t)r^\Delta(t))}dh\right\}r^\Delta(t)\nonumber,
\end{align} we obtain the desired result, except at $t=b$ in the
case that $\rho(b)<b$. To handle this case, we just need to
integrate the last inequality in (\ref{seila8}) from $a$ to $b$ and use \cite[Theorem~1.77, item (ix)]{livro}.
\end{proof}

\begin{theorem}
\label{thm1} Let $u(t)$ and $f(t)$ be nonnegative rd-continuous
functions in the time scales interval
$\mathbb{T}_\ast:=[a,b]_{\mathbb{T}}$ and $\mathbb{T}_\ast^{\kappa}$,
respectively. Let $k(t,s)$ be defined as in \cite[Theorem~1.117]{livro}
in such a way that $k(t,s)$ and $k^{\Delta_1}(t,s)$ are nonnegative for
every $t,s\in\mathbb{T}_\ast$ with $s\leq t$ for which they are
defined (it is assumed that $k$ is not identically zero on
$\mathbb{T}_\ast^{\kappa}\times\mathbb{T}_\ast^{\kappa^2}$). Let $\Phi\in
C(\mathbb{R}^+_0,\mathbb{R}^+_0)$ be a nondecreasing, subadditive
and submultiplicative function, such that $\Phi(u)>0$ for $u>0$ and let $W\in C(\mathbb{R}^+_0,\mathbb{R}^+_0)$ be a nondecreasing
function such that for $u> 0$ we have $W(u)>0$. Assume that $a(t)$
is a positive rd-continuous function and nondecreasing for
$t\in\mathbb{T}_\ast$. If
\begin{equation}
\label{eq0} u(t)\leq a(t)+\int_a^t f(s)u(s)\Delta
s+\int_a^tf(s)W\left(\int_a^s
k(s,\tau)\Phi(u(\tau))\Delta\tau\right)\Delta s,
\end{equation}
for $a\leq\tau\leq s\leq t\leq b$, $\tau, s, t\in\mathbb{T}_\ast$,
then for all $t\in\mathbb{T}_\ast$ satisfying
\begin{equation*}
\Psi(\zeta)+\int_a^{\rho(t)} k(\rho(t),s)\Phi(p(s))\Phi\left(\int_a^s
f(\tau)\Delta\tau\right)\Delta s\in Dom(\Psi^{-1})
\end{equation*}
we have
\begin{multline}
\label{eq11}
u(t)\leq p(t) a(t) \\
+ p(t) \int_a^t f(s)W\left[\Psi^{-1}
\left(\Psi(\zeta) +\int_a^sk(s,\tau)\Phi(p(\tau))
\Phi\left(\int_a^\tau
f(\theta)\Delta\theta\right)\Delta\tau\right)\right]\Delta s \, ,
\end{multline}
where
\begin{gather}
\label{p} p(t) = 1+\int_a^t f(s)e_{f}(t,\sigma(s))\Delta s \, ,\\
\zeta = \int_a^{\rho(b)} k(\rho(b),s)\Phi(p(s)a(s))\Delta s \, , \nonumber \\
\label{p1}\Psi(x) = \int_{x_0}^x\frac{1}{\Phi(W(s))}ds,\ x>0,\ x_0>0 \, ,
\end{gather}
and $\Psi^{-1}$ is the inverse of $\Psi$.
\end{theorem}
\begin{remark}
We are interested to study the situation when
$k$ is not identically zero on
$\mathbb{T}_\ast^{\kappa}\times\mathbb{T}_\ast^{\kappa^2}$.
That comprise the new cases, not considered previously in the literature. The case $k(t,s) \equiv 0$ is studied in
\cite[Th.~3.1]{Pach} and is not discussed here.
\end{remark}
\begin{proof}
Define the function $z(t)$ in $\mathbb{T}_\ast$ by
\begin{equation}
\label{eq1} z(t)=a(t)+\int_a^tf(s)W\left(\int_a^s
k(s,\tau)\Phi(u(\tau))\Delta\tau\right)\Delta s \, .
\end{equation}
Then, (\ref{eq0}) can be restated as
\begin{equation*}
u(t)\leq z(t)+\int_a^t f(s)u(s)\Delta s.
\end{equation*}
Clearly, $z(t)$ is rd-continuous in $t\in\mathbb{T}_\ast$. Using
Gronwall's inequality \cite[Theorem~2.7]{Pach}, we get
\begin{equation*}
u(t)\leq z(t)+\int_{a}^t f(s)z(s)e_f(t,\sigma(s))\Delta s \, .
\end{equation*}
Moreover, it is easy to see that $z(t)$ is nondecreasing in
$t\in\mathbb{T}_\ast$. We get
\begin{equation}
\label{eq2} u(t)\leq z(t)p(t),
\end{equation}
where $p(t)$ is defined by (\ref{p}). Define
$$v(t)=\int_a^t k(t,s)\Phi(u(s))\Delta s,\ t\in\mathbb{T}_\ast^{\kappa}.$$
From (\ref{eq2}), and taking into account the properties of $\Phi$, we observe that
\begin{equation*}
\begin{split}
v(t)&\leq\int_a^t
k(t,s)\Phi\left[p(s)\left(a(s)+\int_a^sf(\tau)W(v(\tau))
\Delta\tau\right)\right]\Delta s\\
&\leq\int_a^t k(t,s)\Phi(p(s)a(s))\Delta s+\int_a^t
k(t,s)\Phi\left(p(s)\int_a^sf(\tau)W(v(\tau))\Delta\tau\right)\Delta
s\\
&\leq\int_a^{\rho(b)} k(\rho(b),s)\Phi(p(s)a(s))\Delta s\\
& \qquad +\int_a^t
k(t,s)\Phi\left(p(s)\int_a^sf(\tau)\Delta\tau\right)
\Phi(W(v(s)))\Delta
s\\
&=\zeta+\int_a^t
k(t,s)\Phi\left(p(s)\int_a^sf(\tau)\Delta\tau\right)
\Phi(W(v(s)))\Delta s \, .
\end{split}
\end{equation*}
Define function $r(t)$ on $\mathbb{T}_\ast^{\kappa}$ by
$$r(t) = \zeta+\int_a^t
k(t,s)\Phi\left(p(s)\int_a^sf(\tau)\Delta\tau\right)
\Phi(W(v(s)))\Delta s \, .$$ Since $p$ and $a$ are positive functions, we have that $\Phi\left(a(s)p(s)\right)>0$ for all $s\in\mathbb{T}_\ast$. Since
$k^{\Delta_1} \geq 0$, we must have $\zeta>0$, hence
$r(t)$ is a positive function on $\mathbb{T}_\ast^{\kappa}$.
In addition, $r(t)$ is delta differentiable on $\mathbb{T}_\ast^{\kappa^2}$ with
\begin{align}
r^\Delta(t)&=k(\sigma(t),t)\Phi\left(p(t)\int_a^t
f(\tau)\Delta\tau\right) \Phi(W(v(t)))\nonumber\\
&\ \ \ +\int_a^t
k^{\Delta_1}(t,s)\Phi\left(p(s)\int_a^sf(\tau)\Delta\tau\right)
\Phi(W(v(s)))\Delta s\nonumber\\
& \label{eq4}\\
&\leq\Phi(W(r(t)))\Biggl[k(\sigma(t),t)\Phi\left(p(t)\int_a^t
f(\tau)\Delta\tau\right)\nonumber\\
&\ \ \ +\int_a^t
k^{\Delta_1}(t,s)\Phi\left(p(s)\int_a^sf(\tau)\Delta\tau\right)\Delta
s\Biggr].\nonumber
\end{align}
Dividing both sides of inequality (\ref{eq4}) by $\Phi(W(r(t)))$, we obtain
\begin{equation*}
\frac{r^\Delta(t)}{\Phi(W(r(t)))}\leq\left[\int_a^t
k(t,s)\Phi\left(p(s)\int_a^sf(\tau)\Delta\tau\right)\Delta
s\right]^\Delta.
\end{equation*}
Let us consider the function $\Psi$ defined by (\ref{p1}).
Delta-integrating this last inequality from $a$ to $t$ and
using Lemma~\ref{lemimp}, we obtain
$$\Psi(r(t))\leq\Psi(r(a))+\int_a^t
k(t,s)\Phi\left(p(s)\int_a^sf(\tau)\Delta\tau\right)\Delta s,$$
from which it follows that
\begin{equation}
\label{seila5} r(t)\leq\Psi^{-1}\left(\Psi(\zeta)+\int_a^t
k(t,s)\Phi(p(s))\Phi\left(\int_a^s f(\tau)\Delta\tau\right)\Delta s\right),\
t\in\mathbb{T}_\ast^{\kappa}.
\end{equation}
Combining (\ref{seila5}), (\ref{eq2}) and (\ref{eq1}), we obtain
the desired inequality (\ref{eq11}).
\end{proof}
If we let $\mathbb{T}=\mathbb{R}$ in Theorem~\ref{thm1}, we get
\cite[Th.~2.1]{motivacao}. If in turn we consider
$\mathbb{T}=\mathbb{Z}$, then we obtain the following result:
\begin{corollary}
Let $u(t)$ and $f(t)$ be nonnegative functions in the time scales
interval $\mathbb{T}_\ast:=[a,b]_{\mathbb{Z}}$ and
$[a,b-1]_{\mathbb{Z}}$, respectively. Let $k(t,s)$ be defined as
in \cite[Theorem~1.117]{livro} in such a way that $k(t,s)$ and
$k^{\Delta_1}(t,s)=k(\sigma(t),s)-k(t,s)$ are nonnegative for
every $t,s\in\mathbb{T}_\ast$ with $s\leq t$ for which they are
defined (it is assumed that $k$ is not identically zero on
$[a,b-1]_{\mathbb{T}_\ast}\times[a,b-2]_{\mathbb{T}_\ast}$). Let
$\Phi\in C(\mathbb{R}^+_0,\mathbb{R}^+_0)$ be a nondecreasing,
subadditive and submultiplicative function such that $\Phi(u)>0$
for $u>0$ and let $W\in C(\mathbb{R}^+_0,\mathbb{R}^+_0)$ be a
nondecreasing function such that for $u> 0$ we have $W(u)>0$.
Assume that $a(t)$ is a positive and nondecreasing  function for
$t\in\mathbb{T}_\ast$. If
\begin{equation*}
u(t)\leq a(t)+\sum_{s=a}^{t-1}
f(s)u(s)+\sum_{s=a}^{t-1}f(s)W\left(\sum_{\tau=a}^{s-1}
k(s,\tau)\Phi(u(\tau))\right),
\end{equation*}
for $a\leq\tau\leq s\leq t\leq b$, $\tau, s, t\in\mathbb{T}_\ast$,
then for all $t\in\mathbb{T}_\ast$ satisfying
\begin{equation*}
\Psi(\zeta)+\sum_{s=a}^{t-2}
k(t-1,s)\Phi(p(s))\Phi\left(\sum_{\tau=a}^{s-1} f(\tau)\right)\in
Dom(\Psi^{-1})
\end{equation*}
we have
\begin{equation*}
u(t)\leq p(t)\left\{a(t)+\sum_{s=a}^{t-1} f(s)W\left[\Psi^{-1}
\left(\Psi(\zeta)+\sum_{\tau=a}^{s-1}
k(s,\tau)\Phi(p(\tau))\Phi\left(\sum_{\theta=a}^{\tau-1}
f(\theta)\right)\right)\right]\right\},
\end{equation*}
where
\begin{eqnarray*}
p(t)=1+\sum_{s=a}^{t-1} f(s)e_{f}(t,s+1) \, ,\\
\zeta=\sum_{s=a}^{b-1} k(b-1,s)\Phi(p(s)a(s)) \, ,\\
\Psi(x)=\int_{x_0}^x\frac{1}{\Phi(W(s))}ds,\ x>0,\ x_0>0 \, ,
\end{eqnarray*}
and $\Psi^{-1}$ is the inverse of $\Psi$.
\end{corollary}

For the particular case $\mathbb{T}=\mathbb{R}$,
Theorem~\ref{thm2} generalizes the result obtained by Oguntuase in \cite[Th.~2.3, 2.9]{og}.

\begin{theorem}
\label{thm2} Suppose that $u(t)$ is a nonnegative rd-continuous
function in the time scales interval
$\mathbb{T}_\ast=[a,b]_\mathbb{T}$ and that $h(t)$, $f(t)$ are
nonnegative rd-continuous functions in the time scales interval
$\mathbb{T}_\ast^{\kappa}$. Assume that $b(t)$ is a nonnegative
rd-continuous function and not identically zero on
$\mathbb{T}_\ast^{\kappa^2}$. Let $\Phi(u)$, $W(u)$,
and $a(t)$ be as defined in Theorem~\ref{thm1}. If
\begin{equation*}
u(t)\leq a(t)+\int_a^t f(s)u(s)\Delta
s+\int_a^tf(s)h(s)W\left(\int_a^s
b(\tau)\Phi(u(\tau))\Delta\tau\right)\Delta s,
\end{equation*}
for $a\leq\tau\leq s\leq t\leq b$, $\tau, s, t\in\mathbb{T}_\ast$,
then for all $t\in\mathbb{T}_\ast$ satisfying
$$
\Psi(\xi)+\int_a^{\rho(t)} b(\tau)\Phi(p(\tau))\Phi\left(\int_a^\tau
f(\theta)h(\theta)\Delta\theta\right)\Delta\tau\in Dom(\Psi^{-1})
$$
we have
\begin{multline*}
u(t)\leq p(t) a(t)\\
+p(t) \int_a^t f(s)h(s)W\left[\Psi^{-1}
\left(\Psi(\xi)+\int_a^sb(\tau)\Phi(p(\tau))\Phi\left(\int_a^\tau
f(\theta)h(\theta)\Delta\theta\right)\Delta\tau\right)\right]\Delta s \, ,
\end{multline*}
where $p(t)$ is defined by (\ref{p}), $\Psi$ is defined by
(\ref{p1}), and
\begin{eqnarray*}
\xi=\int_a^{\rho(b)} b(s)\Phi\left(p(s)a(s)\right)\Delta s.
\end{eqnarray*}
\end{theorem}

\begin{proof}
similar to the proof of Theorem~\ref{thm1}.
\end{proof}

For the remaining of this section, we use the following class of
$S$ functions.

\begin{definition}[($S$ function)~~]
A nondecreasing continuous function
$g:\mathbb{R}^+_0\rightarrow\mathbb{R}^+_0$ is said to belong to
class $S$ if it satisfies the following conditions:
\begin{enumerate}
    \item $g(x)$ is positive for $x>
    0$;
    \item $(1/z)g(x)\leq g(x/z)$ for $x\geq 0$ and $z\geq 1$.
\end{enumerate}
\end{definition}

\begin{remark}
For a brief discussion about this class of $S$ functions, the
reader is invited to consult \cite[Sect.~4]{Beesack}.
\end{remark}

\begin{theorem}
\label{thm:nr:3.6} Let $u(t)$, $f(t)$, $k(t,s)$, $\Phi$ and $W$ be
as defined in Theorem~\ref{thm1} and assume that $g\in S$. Suppose
that $a(t)$ is a positive, rd-continuous and nondecreasing
function. If
\begin{equation}
\label{eq8} u(t)\leq a(t)+\int_a^t f(s)g(u(s))\Delta
s+\int_a^tf(s)W\left(\int_a^s
k(s,\tau)\Phi(u(\tau))\Delta\tau\right)\Delta s,
\end{equation}
for $a\leq\tau\leq s\leq t\leq b$, $\tau, s, t\in\mathbb{T}_\ast$,
then for all $t\in\mathbb{T}_\ast$ satisfying
$$G(1)+\int_a^t f(\tau) \Delta\tau\in Dom(G^{-1})$$
and
$$\Psi(\bar{\zeta})+\int_a^{\rho(t)} k(\rho(t),\tau)\Phi(q(\tau))\Phi\left(\int_a^\tau
f(\theta)\Delta\theta\right)\Delta\tau
\in Dom(\Psi^{-1}),$$ we have
\begin{multline*}
u(t)\leq q(t) \max\{a(t),1\}\\
+ q(t) \int_a^t f(s)W\left[\Psi^{-1}
\left(\Psi(\bar{\zeta})+\int_a^s
k(s,\tau)\Phi(q(\tau))\Phi\left(\int_a^\tau
f(\theta)\Delta\theta\right)\Delta\tau\right)\right]\Delta s \, ,
\end{multline*}
where $\Psi$ is defined by (\ref{p1}),
\begin{eqnarray}
G(x)=\int_{\delta}^{x}\frac{ds}{g(s)},\ x>0,\ \delta >0, \nonumber\\
\label{seila4} q(t)=G^{-1}\left(G(1)+\int_a^t f(\tau) \Delta\tau\right),\\
\bar{\zeta}=\int_a^{\rho(b)}
k(\rho(b),s)\Phi\left(q(s)\max\{a(s),1\}\right)
\Delta s, \nonumber
\end{eqnarray}
and $G^{-1}$ is the inverse function of $G$.
\end{theorem}
\begin{proof}
Define the function
$$z(t)=\max\{a(t),1\}+\int_a^tf(s)W\left(\int_a^s
k(s,\tau)\Phi(u(\tau))\Delta\tau\right)\Delta s.$$ Then, from
(\ref{eq8}) we have that
\begin{equation*}
u(t)\leq z(t)+\int_a^t f(s)g(u(s))\Delta s.
\end{equation*}
Clearly, $z(t)\geq 1$ is rd-continuous and nondecreasing. Since
$g\in S$, we have
\begin{equation*}
\frac{u(t)}{z(t)}\leq 1+\int_a^t
f(s)g\left(\frac{u(s)}{z(s)}\right)\Delta s,
\end{equation*}
or
\begin{equation}
\label{eq:li} x(t)\leq 1+\int_a^t f(s)g(x(s))\Delta s,
\end{equation}
with $x(t)=u(t)/z(t)$. If we define $v(t)$ as the right hand side
of inequality (\ref{eq:li}), we have that $v(a)=1$,
$$v^\Delta(t)=f(t)g(x(t)),$$
and since $g$ is nondecreasing,
$$v^\Delta(t)\leq f(t)g(v(t)),$$
or
\begin{equation}
\label{eq9} \frac{v^\Delta(t)}{g(v(t))}\leq f(t).
\end{equation}
Being the case that $v^\Delta(t)\geq 0$, $\Delta$-integrating
(\ref{eq9}) from $a$ to $t$ and applying Lemma~\ref{lemimp}, we
obtain
$$G(v(t))\leq G(1)+\int_a^t f(\tau)\Delta\tau,$$
which implies that
$$v(t)\leq G^{-1}\left(G(1)+\int_a^t f(\tau)\Delta\tau\right).$$
We have just proved that $x(t)\leq q(t)$, which is equivalent to
\begin{equation*}
u(t)\leq q(t)z(t).
\end{equation*}
Following the same arguments as in the proof of
Theorem~\ref{thm1}, we obtain the desired inequality.
\end{proof}

If we consider the time scale
$\mathbb{T}=h\mathbb{Z}=\{hk:k\in\mathbb{Z}\}$, where $h>0$, then
we obtain the following result.

\begin{corollary}
Let $a,b\in h\mathbb{Z}$,  $h>0$. Suppose that $u(t)$, $f(t)$,
$k(t,s)$, $\Phi$ and $W$ are as defined in Theorem~\ref{thm1} and assume that $g\in S$. Suppose that $a(t)$ is a positive and
nondecreasing function. If
\begin{equation*}
u(t)\leq a(t)+\sum_{s\in[a,t)_{\mathbb{T}_\ast}}
f(s)g(u(s))h+\sum_{s\in[a,t)_{\mathbb{T}_\ast}}
f(s)W\left(\sum_{\tau\in[a,s)_{\mathbb{T}_\ast}} k(s,\tau
)\Phi(u(\tau ))h\right)h,
\end{equation*}
for $a\leq\tau\leq s\leq t\leq b$, $\tau, s,
t\in\mathbb{T}_\ast=[a,b]_{h\mathbb{Z}}$, then for all
$t\in\mathbb{T}_\ast$ satisfying
$$G(1)+\sum_{\tau\in[a,t)_{\mathbb{T}_\ast}} f(\tau)h\in Dom(G^{-1})$$
and
$$\Psi(\bar{\zeta})+\sum_{\tau\in[a,t-h)_{\mathbb{T}_\ast}} k(t-h,\tau )\Phi\left(q(\tau )\right)\Phi\left(\sum_{\theta\in[a,\tau)_{\mathbb{T}_\ast}}
f(\theta)h\right)h\in Dom(\Psi^{-1}),$$ we have
\begin{multline*}
u(t)\leq q(t)\Biggl\{\max\{a(t),1\}\\
+\sum_{s\in[a,t)_{\mathbb{T}_\ast}} f(s)W\left[\Psi^{-1}
\left(\Psi(\bar{\zeta})+\sum_{\tau\in[a,s)_{\mathbb{T}_\ast}}
k(s,\tau )\Phi\left(q(\tau
)\right)\Phi\left(\sum_{\theta\in[a,\tau)_{\mathbb{T}_\ast}} f(\theta
)h\right)h\right)\right]h\Biggr\},
\end{multline*}
where $\Psi$ is defined by (\ref{p1}),
\begin{eqnarray*}
G(x)=\int_{\delta}^{x}\frac{ds}{g(s)},\ x>0,\ \delta >0,\\
q(t)=G^{-1}\left(G(1)+\sum_{\tau\in[a,t)_{\mathbb{T}_\ast}} f(\tau)h\right),\\
\bar{\zeta}=\sum_{s\in[a,b-h){\mathbb{T}_\ast}}
k(b-h,s)\Phi\left(q(s)\max\{a(s),1\}\right)h,
\end{eqnarray*}
and $G^{-1}$ is the inverse function of $G$.
\end{corollary}

\begin{theorem}
Let $u(t)$, $f(t)$, $b(t)$, $h(t)$, $\Phi$ and $W$ be as defined
in Theorem~\ref{thm2} and assume that $g\in S$. Suppose that
$a(t)$ is a positive, rd-continuous and nondecreasing function. If
\begin{equation*}
u(t)\leq a(t)+\int_a^t f(s)g(u(s))\Delta
s+\int_a^tf(s)h(s)W\left(\int_a^s
b(\tau)\Phi(u(\tau))\Delta\tau\right)\Delta s,
\end{equation*}
for $a\leq\tau\leq s\leq t\leq b$, $\tau, s, t\in\mathbb{T}_\ast$,
then for all $t\in\mathbb{T}_\ast$ satisfying
$$\Psi(\bar{\xi})+\int_a^{\rho(t)} b(\tau)\Phi(q(\tau))\Phi\left(\int_a^\tau
f(\theta)h(\theta)\Delta\theta\right)\Delta\tau\in Dom(\Psi^{-1}),$$ we
have
\begin{multline*}
u(t)\leq q(t) \max\{a(t),1\}\\
+ q(t) \int_a^t f(s)h(s)W\left[\Psi^{-1}
\left(\Psi(\bar{\xi})+\int_a^sb(\tau)\Phi(q(\tau))\Phi\left(\int_a^\tau
f(\theta)h(\theta)\Delta\theta\right)\Delta\tau\right)\right]\Delta s,
\end{multline*}
where $\Psi$ is defined by (\ref{p1}), $q(t)$ is
defined by (\ref{seila4}) and
$$\bar{\xi}=\int_a^{\rho(b)} b(s)\Phi\left(q(s)\max\{a(s),1\}\right)\Delta s,$$
\end{theorem}

\begin{proof}
similar to the proof of Theorem~\ref{thm:nr:3.6}.
\end{proof}


\section{An Application}
\label{aplic}

In this section we use Theorem~\ref{thm2} to the qualitative
analysis of a nonlinear dynamic equation.
Let $a,b\in\mathbb{T}$ and consider the initial value
problem
\begin{equation}
u^\Delta(t)=F\left(t,u(t),\int_a^t K(t,u(s))\Delta s\right)\, , \quad t\in\mathbb{T}_\ast^{\kappa}\, , \quad u(a)=u_a \, ,\label{ivp}
\end{equation}
where $\mathbb{T}_\ast=[a,b]_{\mathbb{T}}$, $u\in$
C$_{\textrm{rd}}^1[\mathbb{T}_\ast]$, $F\in$
C$_{\textrm{rd}}[\mathbb{T}_\ast\times\mathbb{R}\times\mathbb{R},\mathbb{R}]$
and $K\in$
C$_{\textrm{rd}}[\mathbb{T}_\ast\times\mathbb{R},\mathbb{R}]$.

In what follows, we shall assume that the IVP (\ref{ivp}) has a
unique solution, which we denote by $u_\ast(t)$.

\begin{theorem}
Assume that the functions $F$ and $K$ in (\ref{ivp}) satisfy the
conditions
\begin{align}
|K(t,u)|&\leq h(t)\Phi(|u|)\label{cond1},\\
|F(t,u,v)|&\leq |u|+|v|\label{cond2},
\end{align}
where $h$ and $\Phi$ are as defined in Theorem~\ref{thm2}. Then,
for $t\in\mathbb{T}_\ast$ such that
$$\Psi(\xi)+\int_a^{\rho(t)} \Phi(p(\tau))\Phi\left(\int_a^\tau
h(\theta)\Delta\theta\right)\Delta\tau\in Dom(\Psi^{-1}),$$ we have the
estimate
\begin{equation}
\label{seila3} |u_\ast(t)|\leq p(t)\left\{|u_a|+\int_a^t
h(s)\Psi^{-1} \left(\Psi(\xi)+\int_a^s
\Phi(p(\tau))\Phi\left(\int_a^\tau
h(\theta)\Delta\theta\right)\Delta\tau\right)\Delta s\right\},
\end{equation}
where
\begin{eqnarray*}
p(t)=1+\int_a^t e_{1}(t,\sigma(s))\Delta s,\\
\xi=\int_a^{\rho(b)} \Phi(p(s)|u_a|)\Delta s,\\
\Psi(x)=\int_{x_0}^x\frac{1}{\Phi(s)}ds,\ x>0,\ x_0>0 \, .
\end{eqnarray*}
\end{theorem}
\begin{proof}
Let $u_\ast(t)$ be the solution of the IVP (\ref{ivp}). Then,
we have
\begin{equation}
\label{seila1} u_\ast(t)=u_a+\int_a^t F\left(s,u_\ast(s),\int_a^s
K\left(s,u_\ast(\tau)\right)\Delta\tau)\right)\Delta s.
\end{equation}
Using (\ref{cond1}) and (\ref{cond2}) in (\ref{seila1}), we have
\begin{align}
|u_\ast|&\leq|u_a|+\int_a^t
\left(|u_\ast(s)|+\int_a^s|K(s,u_\ast(\tau))|\Delta\tau\right)\Delta
s\nonumber\\
&\leq|u_a|+\int_a^t
\left(|u_\ast(s)|+h(s)\int_a^s\Phi\left(|u_\ast(\tau)|\right)\Delta\tau\right)\Delta
s\label{seila2}.
\end{align}
A suitable application of Theorem~\ref{thm2} to (\ref{seila2}), with $a(t)=|u_a|$, $f(t)=b(t)=1$ and $W(u)=u$, yields (\ref{seila3}).
\end{proof}



\end{document}